%
\documentstyle{amsppt}
\magnification \magstep1
\parskip 11pt
\parindent .3in
\pagewidth{5.2in}
\pageheight{7.2 in}
 
%
%
\def \bl{\vskip 11pt}

\def \ni{\noindent}
\def \P{\bold{P}}

\def \C{\bold{C}}
\def \O{\Cal{O}}

\def \lra{\longrightarrow}
%

\centerline {\bf Multiplier Ideals, Vanishing Theorem}
\centerline {\bf and Applications}
\bl
\bl
\bl
\centerline {\smc Lawrence Ein
\footnote {Partially supported by NSF Grant DMS 93-02512}}
 
\ni{\bf \S 0. Introduction.}
 
The purpose of this note is to give a survey of the algebraic
properties of multiplier ideals, and illustrate some of their
applications to classical projective geometry.

Over the past ten or fifteen years, there have been two parallel
developments of the classical Kodaira vanishing theorem. On
the algebraic side, the Kawamata-Viehweg vanishing theorem for
{\bf Q}-divisors has found many important applications, generally
following the model of Kawamata's arguments in [KMM]. From
a more analytic viewpoint, the Nadel vanishing theorem [N]
involving sigular metrics has likewise proven of great utility.
While a priori the analytic notion is more general, in actual
applications these two approaches are essentially just different
ways of packaging the same information. See Demailly's CIME notes
([De3]) and Koll\'ar's lectures in this volume ([Kol4]) for proofs
of this equivalence. While algebraic geometers are probably more
comfortable with ideas occuring in Kawamata-Viehweg-based arguements,
these arguments tend to be complicated by the necessity of passing
to various resolutions of the ambient
variety. By contrast, the Nadel-based
approach essentially works directly on the variety of interest, which
to some authors may seem more geometrically natural.
One can suspect that the conceptual simplicity so achieved was
one of the factors making possible
the breakthrough of Anghern and Siu [AS]. It certainly was
important to us in studying singularities of theta divisors [EL].
However already in the work of Esnault-Viehweg ([EV1]), it was
realized that one could develop the multiplier-ideal viewpoint
algebraically. The purpose of these notes is to outline this
development, and indicate some applications. In term of actual
content, much of the present material is covered
in much greater generality
from the Kawamta-Viehweg viewpoint in Koll\'ar's lectures ([Kol4]).
Nonetheless, we hope that some readers may profit from seeing
the material presented here in the language of multiplier ideals.
 
We start in $\S 1$ and $\S 2$ by recalling the
basic constructions of multiplier
ideals. Then we study the behaviour of these ideals under various
standard geometric operations. In $\S 3$, we apply these results
to study the singularities of theta divisors. In particular,
we show that an
irreducible theta divisor has only rational
singulariteies. All the results in this
section are joint work with R. Lazarsfeld. In $\S 4$, we study
adjoint linear systems. These results generalize the classical
theorems of Kodaira, Bombieri and Reider ([B] and [R])
on linear systems on
surfaces. In the last section, we give a simple proof
of a theorem of Levine on the invariance of plurigenra under
deformations. I would like
to mention that the idea of this proof is due to Siu. Finally, we
discuss a result of Esnault and Viehweg on the zeros of polynomials
([EV2]).
This paper was orignally conceived as being a joint work with
Lazarsfeld and Siu. We
would like to thank them for sharing with us many of their
ideas and many helpful discussions
 
\ni {\bf \S 1. Basic algebraic constructions of multiplier ideals and
adjoint ideals.}
 
To motivate the algebraic
construction of multiplier ideals, it will be useful
to recall the following well known theorems of Esnault and Viehweg
([EV1]  5.1 and 5.13). We consider the following setup. Let $X$ be a
smooth complex projective
variety and let $B$ be a Cartier divisor on $X$. Suppose that
$B \equiv \Delta$, where $\Delta$ is a {\bf Q}-divisor of the form,
$\Delta = \sum a_j F_j$, where the $F_j$'s are distinct smooth
irreducible divisors. We also assume that
the support of $\Delta$ is in normal
crossing, and $0 < a_j < 1$ for all $j$. In other words, $\Delta$
is a boundary divisor.
 
\ni {\bf Theorem 1.1.}
Let $D$ be any effective divisor on $X$
such that the support of $D$ is contained in the support of
$\Delta$. Then the natural map,
$$H^i(\O_X(K_X+B)) \lra H^i(\O_X(K_X+B+D))$$
is injective for every $i$.
 
The next proposition allows us to replace $X$ by another
birational model.
 
\ni {\bf Proposition 1.2.} (a) Let $X$, $B$, and $\Delta$ be as
in Theorem 1.1.
Suppose that $g: Y \lra X$ is a proper birational map such that the
the support of $g^*(\Delta)$ is again a divisor in normal crossing.
Then $g_*(\O_Y(K_Y+g^*B-\lbrack g^*\Delta \rbrack)) = \O_X(K_X+B)$,
and $R^ig_*(\O_Y(K_Y+B- \lbrack \Delta \rbrack)) = 0$
for $i > 0$.
 
\ni (b) Suppose that there is a nef and big divisor
$H$ on $X$, such that
$supp (H) \subset supp (\Delta)$.
Then $H^i (\O_X(K_X+B)) = 0$ for $i > 0$.
 
These theorems are very useful. For instance, they imply the
well known vanishing theorem of Kawamata and Viehweg and they also
give Koll\'ar's
theorem on higher direct images of dualizing sheaf.
 
Now we'll recall the algebraic construction and the basic
properties of multiplier ideals. We consider the following setup:
 
Let $X$ be a smooth varietiy and $G$ be an effective
{\bf Q}-divisor on $X$.
We construct an embedded resolution for $G$,
$$f: Y \lra X.$$
We consider the following {\bf Q}-divisor on $Y$:
$$R = K_Y - f^*(K_X+G) = \sum a_j F_j,$$
where we assume that the
$F_j$'s are distinct irreducible smooth divisors and their
supports are in normal crossing. Then
$$\lceil R \rceil = \sum \lceil a_j \rceil F_j.$$
We write
$$\lceil R \rceil = P-N \ \ \text{or}  \ \ R= P-N-\Delta $$
where $P$ and $N$ are effective integral
divisors with no common components, and
$\Delta$ is an effective {\bf Q}-divisor with
all its coefficients between
0 and 1 and whose suppport is a divisor in normal crossing.
We observe that
$$K_Y+\Delta \equiv f^*(K_X+G)+P-N.$$
We
note that $P$ is $f-$exceptional. By a well known lemma of
Fujita, we know
that $f_*(\O_P(P)) = 0$ ([KMM]). Then
$$ f_*(\O_Y(P-N)) = f_*(\O_Y(-N)) \subset f_*(\O_Y)= \O_X.$$
Hence, $f_*(\O_Y(P-N))$ is an ideal sheaf. We call this the multiplier
ideal of G. Suppose that
$Z(G)$ is the scheme defined by
this ideal. We will denote the mulitplier ideal by $I_{Z(G)}$.
\vskip5pt
 
\ni { \bf Remarks. 1.3.}
 
\ni (a) By some fairly standard methods,
one can check that the ideal
$I_{Z(G)}$ is independent of the
choice of the resoluton. Let $p$ be a point in $X$. We say that
the ideal
$I_{Z(G)}$ is trivial at $p$, if $p$ does not belong to the scheme
$Z(G)$. This
is equivalent to saying
that $f^{-1}(p)$ does not intersect the divisor N.
 
\ni (b) Observe that $F_j \subset N$ if and only if $a_j \le -1$.
This leads to the interpretation that $Z(G)$ puts a scheme structure
on the locus where the pair $(X,G)$ is not log-terminal.
 
\ni (c) If we do not assume that $G$ is an effective {\bf Q}- divisor
in the construction of the multiplier ideal, then $P$ may not be
$f-$exceptional. In this case, $f_*(\O_Y(P-N))$ will give
a fractional ideal. In a similar manner, we say that multipier
fractional ideal is nontrivial at $p$, if $\O_{X,p}$ is not contained
in the localization of this fractional ideal at $p$. Equivalently,
this means that $f^{-1}(p)$ has an nonempty intersection with
the divisor $N$. In this fashion, we construct fractional ideals,
which enjoy many similar cohomological properties as the multipier
ideals.
 
\ni (d) One can in fact carry out the same constrution when
$X$ is normal
and $K_X+G$ is a {\bf Q}-Cartier divisor.
 
The next Proposition follows easily
from the Kawamata-Viehweg-Nadel vanishing theorem or Theorem
1.1 and 1.2.
 
\ni { \bf Propositon 1.4.}
(a) $R^if_*(\O_Y(P-N))=0$ for $i>0$.
 
\ni (b) Assume that $X$ is complete and $A$ is a Cartier divisor
on $X$ such
that $A-(K_X+G)$  is nef and big. Then
$H^i(\O_X(A) \otimes I_{Z(G)})=0$ for $i>0$.
 
\ni Proof. (a) Observe that
$$ P-N = K_Y+ f-(nef \ \ \text{divisor})+\Delta, $$
where $\Delta$ is a boundary divisor
with normal crossing support. Now (a) follows from the
Kawamata-Viehweg vanishing theorem [KMM].
 
\ni (b) It follows from (a) and the Leray spectral sequence that
$$H^i(f^*(\O_Y(A+P-N))) = H^i(\O_X(A) \otimes I_{Z(G)}).$$
Now (b) follows by applying the Kawamata-Viehweg
vanishing theorem on $Y$.\qed
 
Next, we would like to reveiw
the notion of adjoint ideals to study
hypersurface singularities [EL3]. Its construction is very
similar to the construction of
multiplier ideals. Let $X$ be a smooth variety.
Let $H$ be a reduced effective
Cartier divisor in $X$.
Let $f: Y \lra X$ be an embedded resolution and $F$ be the proper
transform of $H$ in $Y$. Then we can write
$$K_Y+F-f^*(K_X+H) = P-N,
\tag{1.5.1}$$
where $P$ and $N$ are effective divisors with no common components.
Note that
every component of $P$ is $f$-exceptional. Then the adjoint ideal
$J = J_H$ is defined to be
$$ J = f_*(\O_Y(P-N)) = f_*(\O_Y(-N)) \subset f_*\O_Y = \O_X.$$
We observe that
$$\O_Y(P-N)|_F = K_F-f^*(K_H).$$
This shows that $J$ is trivial if and only if $H$ has only canonical
singularities. Since $H$ is Gorenstein, $H$ has only canoncial
singularities if and only if it has only rational singularities (cf.
[Kol3, 11.10]).
 
\ni {\bf Proposition 1.5.}
Let $J$ be the adjoint ideal of $H$.  There is the following exact
sequence:
$$0 \lra \O_X(K_X) \lra \O_X(K_X+H) \otimes J \lra f_*(\O_F(K_F))
\lra 0.
\tag{1.5.2}$$
Moreover, $J = \O_X$ if and only if $H$ is normal and it
has only rational singularites.
 
\ni{\bf Proof.} There is the following exact sequence on $Y$:
$$0 \lra \O_Y(K_Y) \lra \O_Y(K_Y+F) \lra \O_F(K_F) \lra 0.$$
We note that $f_*(\O_Y(K_Y))=\O_X(K_X)$ and $R^1f_*(\O_Y(K_Y))=0.$
Also $f_*(\O_Y(K_Y+F)) = \O_X(K_X+H) \otimes J$ by 1.5.1.
We see that 1.5.2 follows from pushing forward the above exact
sequence under $f$. Since the map $f: F \lra H$ factors through
the normalization of $H$, one sees that $f_*\omega_F = \omega_H$
only if $H$ is normal. We have already observed the last statement.
\qed
 
\ni {\bf Remark 1.6.} (a) One may view from (1.5.2) that the
adjoint ideal measures the failure of adjunction.
 
\ni (b) The above construction can also be carried out when $X$ has
only normal Gorenstein canonical singularities.
\vskip 5pt
 
\ni{\bf \S 2. Geometric properties of multiplier ideals. }
 
\vskip 5pt
 
In this section, we will investigate the properties of the multiplier
ideals under simple geometric operations such as taking
hyperplane sections, specialization, and finite ramified coverings.
 
Keeping notations as in \S 1, let $H$ be a smooth
irreducible hypersurface in $X$.
Let $\alpha$ be the coefficent of $H$ in $G$. We will assume that
$0 \le \alpha < 1$.
By Remark 1.3(a),
the construction of the mulitiplier ideal is independent of the choice of
the resolution. Then we may
assume that $f^*H$, the exceptional divisors of $f$, and $f^*(G)$
are all in normal crossing in $Y$.
Let $F$
be the proper transform
of $H$ in $Y$.
We set $G' = G - \alpha H$
and $G'_H= G'|_H$. We observe that $K_H+G'|_H$
is a {\bf Q}-Cartier divisor on $H$.
Then we can consider the multiplier ideal determined by $G'_H$ in $\O_H$.
The following comparison of $I_{Z(G)}$ and $I_{Z(G'_H)}$
is essentially due to Esnault and Viehweg ([EV1], Proposition 7.5).
 
\ni { \bf Prososition 2.1.}
$I_{Z(G'_H)} \subset Im(I_{Z(G)} \lra \O_H).$
 
\ni {\bf Proof.}
Using the notations as before, we write $f^*(H)-F = \sum h_jF_j$.
Let $R' = K_Y-f^*(K_X)- \sum h_jF_j -G' $.  We note that
$R'|_F = K_F - f^*(K_H-G'_H)$. Then
$f_*(\O_F(\lceil R'|_F \rceil) = I_{Z(G'_H)}$. Next we
consider the divisor,
$$R''= R'-F = K_Y-f^*(K_X-G'-H).$$
Since $\alpha < 1$,
we see that
$$\lceil R'\rceil \subset \lceil R \rceil.
\tag{2.1.1}$$
Consider the exact sequece,
$$0 \lra \O_Y (\lceil R'' \rceil) \lra \O_Y(\lceil R' \rceil)
\lra \O_F(\lceil R' \rceil |_F) \lra 0.$$
By Proposition 1.2, $R^1f_*(\O_Y(\lceil R'' \rceil))=0.$
Set $J= f_*(\lceil R' \rceil ).$  Then $J$
maps onto $I_{Z(G_H)}$. By (2.1.1), $J$ is a subsheaf of $I_{Z(G)}$.
This completes the proof of Proposition 2.1.  \qed
\vskip 5pt
 
\ni { \bf Remarks.2.2.}
 
\ni (a) Proposition 2.1 implies that if $I_{Z(G'_H)}$ is trivial at
a point $p$ in $H$, then $I_{Z(G)}$ is also trivial at $p$.
 
\ni (b) In Propostion 2.1, we can make the same conclusion when
$H$ is only assumed to be irreducible noraml and Cartier.
 
The following is a well known criterion to see whether the multiplier
is nontrivial. See [Dem2] and [EV1].
 
\ni {\bf Corollary 2.3.}
Let $X$ be a smooth n-dimensional
variety and $G$ be an effective $Q-$divisor on $X$.
Let $p$ be a point in $X$.
 
\ni (a) If $Mult_p(G) \ge n$, then the multiplier ideal of $G$ is
nontrivial at $p$.
 
\ni (b) If $Mult_p(G) < 1$, then the multiplier ideal of $G$ is trivial
at $p$.
 
\ni Proof. (a) This follows from a
simple calcuation by blowing up $X$ at
$p$.
 
\ni (b) This is clearly true when $n=1$. For $n > 1$,
if $H$ is a general hyperplane
section of $X$ through $p$, then the multiplier of ideal of $G|_H$
is trivial by induction. Hence the multiplier ideal of $G$ is also
trivial at $p$ by 2.2 (a).
\vskip 10pt

\ni { \bf Definition 2.4.} Let $p$ be a point in $X$.
We say that $G$ is critical at $p$, if
 
\ni (a) $p$ is in $Z(G)$ and $p$ is not in
$Z(\lambda G)$, for any $0 < \lambda < 1$.
 
\ni (b) There is a unique irreducible component $F$ of $N$, such that
$N= F+N_1$ where $N_1 \cap f^{-1}(p)$ is empty.
Then $F$ is called the cirtical component of $G$ at $p$
and $f(F)= Z$ is
called the critical variety of $G$ at $p$.

\ni {\bf Remark 2.5}
 
Suppose that $G$ is an effective {\bf Q}-divisor in satisfies 2.4 (a).
In genereal, $G$ may not be critical at $p$. We would like to
sketch a well known arguement showing that we can perturb the divisor
$G$ a little bit to obtain a new divisor $G'$ which is critical at $p$.
 
We'll assume that $G$ is ample and $f$ is a projective morphism.
We can find an ample {\bf Q}-divisor in $Y$ of the following form:
$$A = f^*(\alpha G) - \sum \beta_j F_j,$$
 
\ni where $0 < \beta_j << \alpha << 1$. By Bertini's theorem,
we can find a sufficently
large and divisible integer $m$ with the following properties.
 
\ni (a) There is a smooth irreducible divisor $H \equiv mA$ and
$H + \sum F_j$ is a divisor of normal crossing.
 
\ni (b) $m\beta_j$ are integers.
 
\ni (c) $m\alpha G$ is an integral divisor.
 
Observe that $|H + m\sum \beta_j F_j| = |f^*(m\alpha G)|$.
So $H + m\sum \beta_jF_j = f^*D$ for some effective divisor $D$
in $X$. Observe that $\frac{1}{m} D \equiv \alpha G$. Now consider
the {\bf Q}-divisor $G_1 = G + \frac{1}{m} D$. By choosing $\beta_j$
sufficently generic, we may assume there is a small positive rational
number $\epsilon$, such that $G' = (1- \epsilon)G_1$ is critical at $p$.
 
The following result has also been observed by Kawamata independently
([Ka2]).
 
\ni { \bf Proposition 2.6.}
Suppose that $G$ is critical at $p$. Let
$F_1$ be its critical component and $Z$ be its critical variety.
Then the natural projection map
$$f: F_1 \lra Z,$$
has a trivial Stein factorization near $p$.  Then $Z$
is normal at $p$.In particular, if $Z$ is a curve, then $Z$ is
smooth at $p$.
 
\ni Proof. After replacing $X$ by an open neighbohood of $p$, we
may assume that $\lceil R \rceil = P - F_1$. By a theorem of Fujita,
we know that $f_*(\O_P(P)) = 0$ ([KMM]).
This implies that
$f_*(\O_Y) = f_*(\O_Y(P)) = \O_X.$ By Proposition 1.2,
$f_*(\O_Y(P))$ maps onto $f_*(\O_{F_1}(P)).$ This implies that
the natural map
$\O_Y \lra f_*(\O_{F_1})$ is onto. We conclude that
$f_*(\O_{F_1}) = \O_Z$, and the Stein factorization of the map from
$F_1$ to $Z$ is trivial.          \qed
 
It is an interesting problem
to understand the singularity of $Z$ at $p$.
It is well known that if $Z$ is of codimension $1$, then $Z$ has only
log-terminal singularities. More recently, Kawamata has shown that
if $codim (Z) = 2$, then $Z$ again only has log-terminal singularity
[Ka2]. One may ask that whether this is true in general.
 
The following result was first shown by Angehrn and Siu using analytic
methods [AS]. See also [Kol4].
 
\ni { \bf Proposition 2.7.}
Let $X$ be a smooth variety and let $T$ be a smooth curve.
We consider the product variety $X \times  T$. Let $G$
be an effective {\bf Q}-Cartier divisor  on $X \times  T$.
We assume that the support
of $G$ does not contain any fiber   $X_t =X \times  \{ t \}$,
for $t \in T$.
Denote by $G_t$
the restriction of $G$ to $ X_t$. Let
$$ s: T \lra X \times  T $$
be a section.
Assume that the mulitplier ideal of $G_t$
is nontrivial at the point $s(t)$
for a general $t$. If $X_0 = X \times  \{ 0 \}$ is the special fiber, then
the multiplier ideal of $G_0$ is also nontrivial at $s(0)$.
 
\ni Proof. We construct an embedded resolution
$f: Y \lra X \times T$. By the theorem of generic smoothness, we see
that for generic $t$ in $T$, $f^{-1} (X_t) \lra X_t$ is also an
embedded resolution. This implies that $Z(G) \cap X_t = Z(G_t)$ for
generic $t$. If $s(0)$ is not in $Z(G_0)$, then $s(0)$ is
not in $Z(G)$ by Proposition 2.1.
This means that $s(t)$ is not in $Z(G)$ for generic $t$.
This is a contradiction. \qed
 
\vskip 15pt
 
The following Proposition discusses the properties of a multiplier
ideal, when the {\bf Q}-divisor is pulled back by a generic finte
map. See also [Kol4].
 
\ni { \bf Proposition 2.8. }
Let $X$ and $M$ be two smooth irreducible varieties. Let
$G$ be an effective
{\bf Q}-divisor on $X$. Let $\phi: M \lra X$ be a proper generically
finite map. Let $q$ be a point in $M$. Then the fractional multiplier
ideal of the {\bf Q}-divisor $\phi^*G - K_{M/X}$ is nontrivial at $q$
if and only if the multiplier ideal of $G$ is nontrivial at $\phi(q)$.
 
\ni {\bf Proof.}
We can construct embedded resolutions, $f: Y \lra X$
and $g: W \lra M$. We may suppose that $\phi$ extends to a proper
gernerically finite map $\phi': W \lra Y$. Let $R= K_Y-f^*(K_X+G)$
and $R_1 = K_W - \phi^{'*}f^*(K_X+G)$. Then the multiplier
ideal of $G$ is given by $f_*(\O_Y(\lceil R \rceil)$. Similarly the
fractional multiplier ideal for $\phi^*G-K_{M/X}$ is given by
$g_*\O_W (\lceil R_1 \rceil)$.
Let $F$ be an irreducible component of
$R$ with coefficent $a$ and let $E$ be an irreduible divisor in $W$
which maps onto $F$. Denote by $m$ the coefficent of $E$ in $K_W/Y$
and $b$ the coefficent of $E$ in $R_1$. Then
$$b= a(m+1)-m.
\tag{2.7.1}
$$
We note
that $a \le -1$ if and only if $b \le -1$. Convesely if $E$ is an
irreducible component of $R_1$, after further blowing up of $Y$ and
$W$, we may assume that $\phi'(E)$ is a divisor in $Y$. It follows
that the multiplier ideal of $G$ is nontrivial at
$\phi(q)$ if and only if the fractional multiplier ideal of
$\phi^*G-K_{M/X}$is nontrivial at $q$.\qed
 
\ni { \bf Example 2.9 }
(Demailly)
Let $D$ be the divisor in $\C^n$ defined by the equation
$x_1^{d_1}+ x_2^{d_2}+...+ x_n^{d_n}.$  Assume that
$\frac{1}{d_1}+\frac{1}{d_2} + ... + \frac{1}{d_n} = \lambda < 1.$
Then $\lambda D$ is critical at the origin. One can check this by
considering the finite map $\phi: \C^n \lra \C^n$ given by
$(x_1^{m_1}, x_2^{m_2}, ..., x_n^{m_n})$, where
$m_i = \frac{d_1d_2...dn}{d_i}$.Then $\phi^*D$ is define
by the Fermat equation $x_1^d+x_2^d+...+x_n^d$ where $d=d_1...d_n$.
The relative canonical divisor of $\phi$ is defined by
$x_1^{m_1-1}...x_n^{m_n-1}$. After blowing up the origin, we
obtain an embedded resolution of $\phi^*D$ and $K_{\phi}$. Now
2.7 follows from 2.6 by a simple calculation. See [Dem3] for an
arguement using the analytic techniques and see
also [Kol4] for
another argument using weighted blowing up.
 
Observe that if $d_i >> d_1$
for $i> 1$, then the multiplicity of $\lambda D$ at the origin is
of the form of $1+\epsilon$ and it is an isolated singularity.
This shows that Proposition 1.5(b) is essentially the
best possible result.
 
\vskip 10pt
\ni{\bf \S3 Singularties of theta divisors}
 
In this section, we'll apply the techniques from the first two
sections to study singularties of theta divisors. The results in
this section are my joint work with R. Lazarsfeld. For more details
see [EL3]. In the following, we let $(A, \Theta)$ be a
principally polarized abelian
variety. First we note the following well known result.
 
\ni{\bf Propostion 3.1.} Let $Z$ be a nonempty proper closed
subscheme of $A$. Then
$$H^0(\O_A(\Theta)\otimes P \otimes I_Z) = 0 \ \
\text{for a general}\ \ P \in Pic^0A.$$
 
\ni{\bf Proof.} We first note that $\O_A(\Theta) \otimes P$ is linearly
equivalent to a translate of $\Theta$. Now the Proposition
follows from $h^0(\O_A(\Theta) \otimes P) =1$
and the fact that
a general translate of the divisor $\Theta$ does not contain $Z$.
\qed
 
The following is an extension of a theorem of Koll\'ar to the
pluritheta divisors, as proposed in [Kol2, Problem 17.15].
 
\ni { \bf Propostion 3.2. }
Let $(A, \Theta)$ be a p.p.a.v., and for a positive integer
$m$ fix any divisor $D \in |m\Theta|$. Then the pair
$(A, \frac{1}{m}D)$
is log-canonical. In particular for each positive integer $k$ the set,
$$\Sigma_{mk}(D) = \{x\in D| \text{mult}_x(D) \ge mk \}$$
is of codimension greater or equal to $k$.
 
\ni Proof. Suppose for contradiction that $(A, \frac{1}{m}D)$
is not log-canonical. Then for sufficently small positive
$\epsilon$, the multiplier ideal $I_Z$ for the divisor
$G = (1-\epsilon) \frac{1}{m}D$ is nontrivial.
By Propositon 1.4,
$$H^i(\O_A(\Theta) \otimes P \otimes I_Z) = 0 \ \ \text{for all} \ \
i > 0\ \  \text{and} \ \ P \in \text{Pic}^0(A).
\tag{3.2.1}
$$
It follows from Proposition 3.1 that
$\chi(\O_A(\Theta) \otimes P \otimes I_Z) =0$
for general $P$. This implies that
$\chi(\O_A(\Theta) \otimes P \otimes I_Z) = 0$ for all $P$.
We conclude that
$$H^i( \O_A(\Theta) \otimes P \otimes I_Z) = 0 \ \ \text{for all}
\ \ i \ge 0 \ \ \text{and for all}  \ \ P \in \text{Pic}^0(A).
\tag{3.2.2}
$$
It follows from Mukai's theory [Muk] of the Fourier functor that
$\O_A(\Theta) \otimes I_Z =0$. Thus we have a contradiction.
\qed
 
\ni{\bf Theorem 3.3.} Assume that $\Theta$ is irreducible. Then
the adjoint ideal of $\Theta$ is trivial. In particular, $\Theta$
is normal and has only rational singularites.
 
\ni{\bf Proof.} Let $f: Y \lra A$ be an embedded resolution.  Let
$W$ be the proper transform of $\Theta$ in $Y$. Let $J$ be the
adjoint ideal of $\Theta$. It follows from 1.5.2 that we have the
following exact sequence:
$$0 \lra P \lra \O_A(\Theta) \otimes P \otimes J \lra
f_*\omega_W \otimes P \lra 0
\tag{3.3.1}
$$
for all $P \in Pic^0(A)$.
By the generic vanishing theorem of
Green and Lazarsfeld [GL],
$$H^i(f_*(\omega_W \otimes P)) =0 \ \ \text{for} \ \ i > 0\ \
\text{and for generic} \ \ P \in
\text{Pic}^0(A).$$
If $J$ is nontrivial, then by Propostion 3.1,
$h^0(\O_A(\Theta) \otimes P \otimes J) = 0$ for generic $P$. We
conclude that $H^0(f_*(\omega_W) \otimes P) = 0$ for generic
$P$. This implies that $\chi(f_*(\omega_W \otimes P) = 0$ for
general $P$. This shows that $\chi(f_*(\omega_W)) = 0$. Since
$\Theta$ is an irreducible ample divisor, it is well known that
$\Theta$ is of general type.
It follows from a theorem of Kawamata and Viehweg that
$\chi(f_*(\omega_W)) > 0$ [KV]. This gives a contradiction. \qed
 
\ni {\bf Remark 3.4.} More generally, Lazarsefeld and I ([EL3])
have shown the following. If $X$ is an irreducible closed
subvariety of an
abelian variety. Let $Y$ be a resolution of singularites of $X$.
Then $X$ is
of general type if and only if $\chi(\omega_Y) > 0$. This result
was
conjectured by Koll\'ar.
\vskip 10pt
 
\ni{\bf \S4 Adjoint linear systems}
 
Let $X$ be a smooth $n-$dimensional complex projective variety
and $L$ be an ample line bundle on $X$. In classical projective
geometry, the adjoint linear system $|K_X+L|$ and the pluricanonical
systems $|mK_X|$ play very important roles in studying
the properties of curves
and surfaces. One would expect that they would be important in
studying higher dimensional varieties as well.
The general philosophy
is that if $L$ is \lq\lq sufficently\rq\rq ample, then the adjoint linear
system $|K_X+L|$ should be free or very ample. More precisely,
we expect that if $|K_X + L|$ is not free or very ample, then
$X$ has a subvariety of low degree with respect to $L$. For instance,
we have the following conjectures of Fujita.
 
\ni (1) $|K_X+(n+1)L|$ is free.
 
\ni (2) $|K_X+(n+2)L|$ is very ample.
 
\ni (3) If $X$ is a smooth minimal $n-$fold of general type,
then $|(n+2)K_X|$ is free and $|(n+3)K_X|$ is very ample on
the canonical model of $X$.
 
One may conjecture the following refinement of (1).
 
\ni (4) Suppose that $L^n > n^n.$Also assume that for every
irrdeucible subvariety $Z$ of $X$,
$$L^{dim Z} \cdot Z \ge (n)^{dim Z}.$$
Then the linear system $|K_X+L|$ is free.
 
In the last few years, these questions have generated a great
deal of work.
In this section, we'll describe some of the recent results on
these problems.
When $X$ is a surface, these conjectures hold by the results of
Bombieri and Reider ([B] and [R]).
In higher dimension, the first breakthrough is
due to Demailly [Dem1].  More specifically,
he shows that $|2K_X+12n^nL|$ is very ample using some very powerful
analytic techniques. Using cohomological methods developed by Kawamata
and Shokurov, the author and Lazarsfeld [EL1]
proved  that the freeness part of
Fujita conjecture is true for threefolds. In a very recent work,
Kawamata [Ka1] has proved a similar result for fourfolds.
Recently, Lee has shown that if $X$ is minimal Gorenstein threefold of
geneal type, then $|5K_X|$ gives a birational morphism [Lee].
Ein and Lazarsfeld and independently Helmke have shown that (4) holds
when $X$ is a threefold [H]. See [Lee] for more precise results for
threefolds. In general for a smooth $n-$dimensional
projective variety $X$,
the work of Angehrn and Siu shows that if
$L^{dim Z} \cdot Z >\binom{n+1}{2}^{dim Z}$ for every subvariety
$Z$ in $X$, then $|K_X+L|$ is free. Furthermore, Koll\'ar has
extended the theorem to singular varieties with only log-terminal singularties.
See [Kol4] for more details.
 
In this section, we would describe how we can
use multiplier ideals to solve some of these problems on adjoint
linear systems. This approach is inspired by the work of Anhern and Siu.
The strategy is fairly simple.
Suppose $p$ is a given point in $X$.
We would like to construct an effective
{\bf Q}-divisor $G$ which is linearly equivalent to $\lambda L$
where $\lambda < 1$, such that the multiplier scheme Z(G) defined by
$G$ is zero dimensional at $p$.
Then the vanishing theorem for multiplier ideals will imply that
the restriction map,
$$H^0(\O_X(K_X+L)) \lra H^0(\O_X(K_X+L)|_Z) $$
is surjective. This in turn implies that we can find a section of
$\O_X(K_X+L)$ that does not vanish at $p$. In order to start,
by the Riemann-Roch theorem, if $L^n > n^n$ , then
we can construct an effective {\bf Q}-divisor $G$ such that $G$ is
equivalent to $\lambda_1L$ with $\lambda_1 < 1$ and $Mult_p(G) \ge n$.
Then the multiplier ideal of $G$ is nontrivial at $p$.
After replacing $G$ by a smaller multiple of $G$ and adding a
small pertubation term as in 2.5, we may assume
that $G$ is critical at $p$.
The difficulty is that
in general we do not know that the critical variety
of $G$ at $p$ is zero-dimensional.
The crucial new idea of Angehrn and Siu
is to construct a new {\bf Q}-divisor
$G'$ of the following form,
$$G' = (1-\epsilon) G +D', $$
such that $G'$ is equivalent to $\lambda' L$, where
$\lambda' < 1$.  Furthermore, $G'$ is critical at $p$ and
the critical variety of $G'$ is a proper subset of the critical
variety of $G$. After repeating this process
a finite number of times, we would then be able to construct
an effective {\bf Q}-divisor with the property that its multiplier
scheme is of zero-dimensional.
In the rest of this section, we'll describe a bit more of the
technical details in constructing the {\bf Q}-divisor $G'$ as above.

First we'll introduce an invariant $def_p(G)$, the
deficit of G. This number roughly speaking is a measure on
the difficulty in constructing the new {\bf Q}-divisor $G'$ as
above with a 0-dimensional multiplier scheme at $p$.
We consider the following setup. Let $X$ be a smooth variety and
$p$ be a point in $X$. Let $G$ be an effective {\bf Q}-divisor in $X$.
First we will assume that the multiplier ideal of $G$ is trivial at $p$.
Let $\pi: X' \lra X$ be the blowup of $X$ at $p$
and $E \subset X'$ be the
exceptional divisor. We construct an embedded resolution for $G$,
$f: Y \lra X$, which factors through $\pi$.
We write the factorization
as $f = \pi \circ g$. Suppose that $f^*(G) = \sum g_j F_j$,
$K_{Y/X} = \sum b_j F_j$,
$K_Y - f^*(K_X+G) = \sum a_j F_j$ and
$g^*(E) = \sum e_j F_j$. Now we define the deficit of $G$ as
$$def_p(G) = \inf_{f(F_j)=p} \{ \frac{a_j+1}{e_j} \}.  $$
 
\ni One checks easily that $def_p (G) \le c$ if and only for every
effective {\bf Q}-divisor $D$ where $Mult_p D \ge c$, we have
$p \in Z(G+D)$. In particular, the definition for the deficit is
independent of the choice of the embedded resolution.
Next we consider the case that when
$G$ has a nontrivial multiplier ideal at $p$, but the multiplier
ideal of $(1-t)G$ is trivial at $p$ for any $t > 0$.
In this case, we define the
deficit of $G$ at $p$ as:
$$ def_p(G) = \lim_{t \to 0^+} def_p((1-t)G).$$
For $t>0$, one notes that
$def_p((1-t)G)-t$Max$(\frac{g_j}{e_j}) \le def_p(G) \le def_p((1-t)G).$
 
The following Proposition follows fairly immediately from the
definition.
 
\ni {\bf Proposition 4.1.} (a) $def_p(G) \le dim(X) - mult_p(G)$.
 
\ni (b) If $G$ is critical at $p$ and if the critical variety is
zero-dimensional at $p$, then $def_p(G) = 0$.
 
In the following, we will assume that $G$ is critical at $p$.
We may achieve this by perturbing the divisor $G$ a bit as in 2.5.
Alternatively, one may work with reducible critical varieties as
in [Ka1].
We will investigate the behaviour of the deficit function when we
restrict $G$ to a general hyperplane section $H$ through the
given point $p$. Let $\Lambda$ be the proper tranform of $H$
in $Y$. Then $f^*(H) = \Lambda + g^*(E)$. By adjunction formula,
$$K_{\Lambda/H} = (K_{Y/X} - g^*(E))|_{\Lambda}.$$
Let $\phi = f|_{\Lambda}$. By Bertini's theorem, we may assume that
$f^*(G)$ intersects $\Lambda$ transversely. If we write
$$\lceil K_{Y/X} - g^*(E) - f^*(G) \rceil = P_1 - N_1,
\tag{4.2.1}
$$
where $P_1$ and $N_1$ are effective divisors with no common
components,
then $\phi_*(\O_{\Lambda}((P_1 - N_1)|_{\Lambda})$ is
the multiplier ideal of
$G|_H$. By Proposition 2.1, we know that the multiplier ideal of
of $G|_H$ is nontrivial.
 
\ni {\bf Proposition 4.2} (a) If $def_p(G) \ge 1$, then
$def_p(G|_H) = def_p(G) - 1$
and $G|_H$ is critical at $p$.
 
\ni (b) Suppose that $Z$ is the critical variety of $G$ at $p$.
Then $def_p(G) \le dim (Z)$
 
\ni {\bf Proof.} (a) One can check this using the formula 4.2.1.
 
\ni (b) If $dim (Z) = 0$, then $def_p(G) =0$ by 4.1. In general,
we induct on the dimension of $Z$. Now (b)
follows from (a)
by considering $G|_H$ and induction.\qed
 
The following  result of Helmke shows that we can use
the deficit to control the
multiplicity of the critical variety [H]. Suppose that $X$ is a
smooth n-dimensional variety
and let $p$ be a given point on $X$. Assume that
$G$ is an {\bf Q}- effective divisor on $X$ and $G$ is critical
at $p$. Let $F \subset Y$ be the critical component of $G$ and
let $Z$ be its critical variety. Suppose that the embedded
dimension of $Z$ at $p$ is $n_1$ and
$dim (Z) = r > 0$. Let $d= def_p(G)$. Set $k= \lbrack r-d \rbrack$.
 
\ni {\bf Theorem 4.3.} (Helmke) $Mult_p(Z) \le \binom{n_1-r+k}{k}$.

\ni {\bf Corollary 4.4} Suppose that the critical variety $Z$
is a surface and
the embedded dimension of $Z$ at $p$ is $n_1$. Then
$Mult_p(Z)\le n_1 - 1$.
 
\ni {\bf Proof.} If $def_p(G)> 1$, then $Z$ is smooth by 2.6 and 4.2.
If we
assume that $def_p (G) \le 1$, then
$Mult_p(Z) \le n_1-1$ by 4.3. \qed

Now we continue our discussion of constructing the divisor
$G$ with a zero dimensional multiplier scheme.
We consider the following definition motivated by the work of Angehrn
and Siu [AS].
 
\ni { \bf Definition 4.5. }
Let $B$ be an effective {\bf Q}-Cartier divisor on the critical variety
$Z$. An effective {\bf Q}-Cartier divisor $D$ on $X$ is said to be a
nice lifting of $B$, if $D$ satisfies the following two properties:
 
\ni (a) $D|_Z = B.$
 
\ni (b) $Z((G + D)|_{X-Z}) = Z(G|_{X-Z}).$
 
\vskip 10pt
 
Let $B$ be an effective {\bf Q}-Cartier divisor on the critical variety
$Z$. Suppose that $D'$ and $D''$ are two nice liftings of $B$.
Let $f: Y \lra X$ be an embedded resolution for $D'$ and $D''$.
Let
$f^*D' = \sum d_j'F_j$, $f^*D''= \sum d_j''F_j$, $f^*G = \sum g_jF_j$,
and $R = K_Y-f^*(K_X+G)= \sum a_jF_j$. After replacing $X$ by an
open neighborhood of $p$, we may assume that
$$\lceil R \rceil = P - F_1$$
where $P$ consists of $f-$exceptional divisors.
Then $f(F_1) = Z$ is a critical variety of $G$.
In applying this result, $G$ usually would be an ample divisor.
As in 2.5, we can perturb $G$ a little bit. So we can assume
the following addition simplifying assumption that
$$\{a_j\} > 0 \ \ \text{for all} \ \ j \ne 1.
\tag{4.6.1}
$$

The following basic lemma will allow us to prove that certain
multiplier ideals are nontrivial. The idea is that the picture is
controlled by the restriction of the divisor to the critical variety.
 
\ni { \bf Lemma 4.6.  } Keep assumptions as above.
 
\ni (a) Suppose that
$p$ is in the multiplier scheme
$Z((1-s)D'+(1-t)G)$ for all sufficently small positive
$s$ and $t$. Then $p$ is also in $Z((1-s')D''+(1-t')G)$ for all
sufficently small positive $s'$ and $t'$. Furthermore
the multiplier scheme
$Z((1-s')D''+(1-t')G)$ is a proper closed subset of the
critical variety $Z$ of $G$.
 
\ni(b) Suppose that $p$ is a smooth point of $Z$ and
$$Mult_p(B) > def_p(G).$$
Then for any nice lifting $D''$ of $B$, the multiplier scheme for
$D''+(1-t)G$ is notrivial at the point $p$
for all sufficently small positive $t$.
 
\ni (c) Let $G'$ be a {\bf Q}-effective divisor of the form
$G' = (1-t)G + D''$ where $t$ is a sufficently small positive number.
Assume that the multiplier ideal of $G'$
is nontrivial at $p$, but the multiplier ideal $(1-s)G'$ is trivial
for all positive $s$. Then
$def_p (G') \le def_p ((1-t)G) - Mult_p (B).$
 
One should think of this lemma as saying that we can use the
divisor $D'$
to study the multiplier ideal
constructed from $D''$ provided only that $D'$ and $D''$ are nice
liftings of the same divisor. For local questions, this allows us to
pretend we are working with a particularly nice lifting.
This Lemma is an algebraic substitute for the computation appearing
in [AS] \S 4.
 
\ni {\bf Proof of 4.6.} Replacing $D'$ by $(1-s)D'$ for sufficently small
and generic $s$, by 4.3.1 we may assume that
$$\{a_j - d_j'\} > 0 \ \ \text{for} \ \ j \ne 1.$$
This means that for sufficently small $t$ and $j \ne 1$,
$$\lceil a_j - d_j'+ t g_j \rceil = \lceil a_j - d_j' \rceil.$$
This implies that
$$\lceil K_Y - f^*(K_X+G+D') \rceil = P' - N' - F_1,$$
and
$$\lceil K_Y - f^*(K_X+D'+(1-t)G') \rceil = P' - N'. $$
Consider the following exact sequence on $Y$:
$$0 \lra \O_Y(P'-N'-F_1) \lra \O_Y(P'-N') \lra
\O_{F_1}(P'-N') \lra 0.$$
By the vanishing theorem, $R^1f_*(\O_Y(P'-N'-F_1) = 0$.
Let $I_{Z'} = f_*(\O_Y(P'-N'))$ be the multiplier ideal of
$D'+(1-t)G$. By our assumtion this is a nontrival ideal at $p$.
It follows that its image in $\O_{Z}$, which is equal to
$f_*(\O_{F_1}(P'-N'))$, is a proper ideal subsheaf at $p$.
In the similar fashion, we write
$$\lceil K_Y-f^*(K_X+G+D'') \rceil = P''-N'' -F_1,$$
and
$$\lceil K_Y-f^*(K_X+D''+(1-t')G) \rceil = P''-N''.$$
Let $I_{Z''} = f_*(\O_Y(P'' - N''))$ be the multiplier ideal
of $D''+(1-t')G$. Since $D'$ and $D''$ both are liftings of $B$,
we observe that $(P'-N')|_{F_1} = (P''- N'')|_{F_1}.$
We cocnlude that the image of the restriction map from $I_{Z''}$
to $\O_{Z}$ is equal to $f_*(\O_{F_1}(P' - N'))$ which is a proper
ideal subsheaf of $\O_{Z}$ at $p$. We conclude that $I_{Z''}$
is nontrivial at $p$. The fact that $Z''$ is a proper closed subscheme
of $Z$ follows from the fact $D'$ is a nice lifting.
 
\ni (b) This is a local question. We will assume that $X$ is affine.
Let $I_{p/X}$ and $I_{p/Z}$ be the maximal ideals correspoinding
to the point $p$ in $X$ and $Z$ respectively. Note the restriction
map $I_{p/X}^k \lra I_{p/Z}^k$ is surjective for each nonnegative
integer $k$. By Bertini's theorem, there is a nice lifting $D'$ of
$B$ such that $Mult_p(D') = Mult_p(B)$. Then
$$Mult_p((1-t)D') > def_p(G)$$
for all sufficently small $t$.
>From the definition of deficit, we see that $p \in Z((1-t)D'+(1-s)G)$
for all sufficently small $s$ and $t$. Now (b) follows from (a).
(c) follows from the proof of (a) and (b).  \qed
 
We'll illustrate how we can apply these techniques to study
adjoint linear systems by giving a quick sketch for the
proof of the theorem of Angehrn and Siu [AS]. See also [Kol4].
 
\ni {\bf Theorem 4.7.} Let $X$ be a smooth n-dimensional complex
projective variety and $L$ be an ample divisor on $X$. Suppose that
for every subvariety $Z$ in $X$,
$$L^{dim Z} \cdot Z > \binom{dim X + 1}{2}^{ dim Z}.$$
Then $|K_X+L|$ is free.
 
First we note the following lemma.
 
\ni {\bf Lemma 4.8} Let $X$ be a smooth complex projective variety
and $p$ be a given point in $X$. Suppose that $G$ is an effective
{\bf Q}-divisor on $X$ such that $G$ is critical at $p$. Let $Z$ be
the critical variety of $G$ at $p$ and $dim (Z) = r$.
Suppose that $A$ is an ample
{\bf Q}-divisor in $X$. We assume that
$$ A^r \cdot Z > r^r.$$
Then there is an effective {\bf Q}-divisor $D$ on $X$ such that
$D$ is equivalent to $A$ and $p \in Z((1-t)G+D)$ for all sufficently
small positive $t$. Furthermore, locally near $p$,
$Z((1-t)G+D)$ is a proper subset $Z$.
 
\ni {\bf Proof.} Let $q$ be a general point of $Z$. Then by the
Riemann-Roch theoerem, we can find a
{\bf Q}-Cartier divisor $B_q$ on $Z$
such that $B_q$ is equivalent to $A|_Z$ and
$$Mult_q(B_q) > r \ge def_q(G).$$
If $D_q$ is nice lifting of $B_q$ which is equivalent to $A$,
then $q \in Z((1-t)G+D_q)$ by 4.6. As we specialize $q$ to $p$,
we obtain the desired {\bf Q}-divisor by Propositin 2.7. \qed
 
\ni {\bf Proof of 4.7.}
The idea of the proof is the following.
We will construct a sequence of {\bf Q}-Cartier divisors
$G_i \equiv \lambda_i L$ where
$\lambda_1 < \lambda_2 ... <1$ satifying the following properties.
 
\ni (a) $G_i$ is critical at $p$ with critical variety $Z_i$.
 
\ni (b) If $dim (Z_i) >0$, then $Z_{i+1}$ is a proper subvariety
of $Z_{i}$.
 
\ni (c) $G_i \equiv \lambda_i L$ and
$\lambda_i < \frac{(n)+ (n-1)+...+(n-i+1)}{m}$, where
$m=\binom{n+1}{2}$.
 
In $k$ steps $(k \le n)$, we will reach the case where $Z_k$ is
zero dimensional at $p$.
 
By the Riemann-Roch theorem, we can find an effective
{\bf Q}-Cartier divisor $G_1$
such that $G_1$ is critical at $p$ with critical variety $Z_1$,
and $G_1$ is equivalent to
$\lambda_1 L$, where
$$\lambda_1 < \frac{n}{m}.$$
If $dim (Z_1) = 0$, then we are done. If $dim (Z_1) = d_1 > 0$, we
continue our construction. Suppose inductively we have constructed
{\bf Q}-Cartier divisor $G_i$ with $G_i$ critical
at $p$ with critcial
subvariety $Z_i$ of dimension $d_i$, where
$$d_i \le n-i.$$
The divisor $G_i$ is equivalent to $\lambda_i L$, where
$$\lambda_i < \frac{(n)+(n-1)+...+(n-i+1)}{m}. $$
If $d_i = 0$, then we are done. Now assume that $d_i > 0$.
Observe that $(\frac{d_i}{m}L)^{d_i}Z_i > d_i^{d_i}$. Using
Lemma 4.8 and adding a small perturbation term,
we find $G_{i+1}$ which is critical at $p$, with
critical subvariety $Z_{i+1}$, which is a proper subvariety of
$Z_i$. Furthermore we can choose $G_{i+1}$ which is equivalent to
$\lambda_{i+1}L$ where $\lambda_{i+1} = \lambda_i + \beta_i$
with $\beta_i < \frac{d_i}{m} \le \frac{n-i}{m}.$
It follows that
$$dim (Z_{i+1}) \le n-i-1.$$
Also,
$$\lambda_{i+1} < \frac{(n)+(n-1)+...+(n-i)}{m}.$$
Now we can complete the proof by induction. \qed
 
\ni {\bf Remark 4.9.} Using Corollary 4.4 and an argument as above,
one can show fairly easily that the freeness part of the Fujita
conjecture is true for threefolds and fourfolds [EL1] and [Ka1].
 
\ni {\bf\S5. Additional applications}
 
In this section, we'll give two further applications. First, we'll
give a simple proof of a result of M. Levine on the invarience
of plurigenra under deformations. The idea of this
proof is due to Siu.
The second
application is to give a simple proof for a result of Esnault and
Viehweg on the relation between singular hypersurfaces and postulation
of a finite set in a complex projective space ([EV2]). We'll need the
following Lemma.
 
\ni {\bf Lemma 5.1} Let $L$ be a big line bundle on a smooth complex
projective variety $X$. Suppose that there is an effective
{\bf Q}-divisor $G$ on $X$ which is linearly equivalent to L,
such that the muliplier ideal of $G$
is trivial. Then $H^i(X, \O_X(K_X+L))= 0$ for $i>0$.
 
\ni {\bf Proof.} Let $m$ be a sufficently large and divisible
positive integer, such that $mG$ is an integral Cartier-divisor.
Let $f: Y \lra X$ be an embedded resolution for the linear system
$|mL|$ and $G$. We can write $f^*(mL) \sim A + B$, such that
$A$ is free
and big and $B$ is an effective divisor in normal crossing which
is the fixed component of $|f^*(mL)|$. Since
$mf^*G\in |f^*(mL)|$, $mf^*G - B$ is effective. By our assumption,
we can write
$$K_{Y/X} - f^*(G) = \sum a_j F_j,$$
where $F_j$'s are distinct irreducible divisor in normal crossing and
$a_j > -1$.
Then $f^*L \sim \frac{A}{m} + \frac{B}{m}$  and
$$K_{Y/X} - \frac{A}{m} - \frac{B}{m} = \sum c_j F_j,$$
where $c_j > -1$.  Then there is an effective $f-$exceptional
divisor $P_1$, such that
$$f^*(K_X+L) + P_1 \sim K_Y+ \Delta_1,$$
where $\Delta_1$ is a boundary divisor with normal crossing support.
By Bertini's theorem, we may assume that
$supp A \subset supp \Delta_1.$ By Propositon 1.2, this implies that
$$H^i(\O_Y(f^*(K_X+L)+P_1)) = 0 \ \ \text{for} \ \ i >0.$$
Now since $f_*(f^*(K_X+L)+P_1) = K_X+L$ and
$R^if_*(f^*(K_X+L)+P)= 0$ for $i > 0$. This implies that
$H^i(\O_X(K_X+L)) = 0$ for $i > 0$. \qed

\ni {\bf Proposition 5.2} (Levine)
Let $f: X \lra T$ be a smooth projective morphism with connected
fibers and $X_0 = f^{-1}(t_0)$ be a closed fiber of $f$. Suppose
that $X_0$ is of general type.
Let $m \ge 2$ be a positive integer. We suppose
that there is a divisor
$D \in |mK_{X_0}|$ and that the multiplier ideal for the
{\bf Q}-divisor $G = \frac{m-1}{m} D \sim (m-1)K_{X_0}$ is trivial.
Then $H^i(\O_{X_0} (aK_{X_0})) = 0$ for $i > 0$ and $2 \le a \le m$.
In particular,
the plurigenra $h^0(\O_{X_t}(aK_{X_t})$ are locally constant for
all $t$ in a neighborhood of $t_0$.
 
\ni {\bf Proof.} Observe that the multiplier ideal of $\lambda G$
is trivial for $\lambda < 1$. We note that
$(a-1)K_{X_0}  \sim \frac{a-1}{m-1} G$. Now the Proposition
follows from 5.1 and the semicontinuity theorem.
\qed
 
\ni{\bf Remark 5.3} (a) Levine actually showed that the plurigenera are
constant without the additional assumption that $X_0$ is of general
type. The observation that one can use multiplier ideals to give
a simple proof of this result, when $X_0$ is of general type, is
due to Siu.
 
\ni (b) The assumption in 5.2 is equivalent to saying
that $G$ is
log-terminal.
 
The final application is a theorem of Esnault and
Viehweg [EV2] on the zeros of the
polynomials. Let $S \subset \P^n$ be a finite subset. Suppose that
there is a hypersurface $D$ of degree $d$ in $\P^n$, such that
$Mult_p D \ge k$ for each $p \in S$.
 
\ni {\bf Proposition 5.4.} There is a hypersurface of degree
$[\frac{nd}{k}]$ that contains $S$.
 
\ni Proof. Let $G = \frac{n}{k}D$. Then $Mult_p G \ge n$ for each
point $p \in S$. Let $I_Z$ be the multiplier ideal of $G$. Then
$S \subset Z$. Let $H$ be the hyperplane class of $\P^n$.
Let $m = [\frac{nd}{k}]$.
Then $(m+1)H - G$ is ample. Since
$K_{\P^n}\sim (-n-1)H$, we conclude by the vanishing theorem that
$H^i (I_Z \otimes \O_{P^n}(t)) =0 $ for $t \ge m - n$ and $i > 0$.
Since the Hilbert polynomial of $I_Z$ is a polynomial in $t$ of degree
less than or equal to $n$, we
can find an integer $t_0$ where $m-n \le t_0 \le m$, such that
$h^0(I_Z \otimes \O_{\P^n} (t_0))$ is nonzero. Since $t_0 \le m$
and $S \subset Z$, this implies the proposition.
\qed
 
\def \ky#1{\item{[{\bf #1}]} }
\parskip 6pt
\parindent .6in
\def \bl{\vskip 14pt}

$$ \text  {\bf  References}$$
\vskip 5pt
 
\ky{AS} U. Angehrn and Y.-T. Siu, {\it Effective freeness and point separation
for adjoint bundles}, Inv. Math. {\bf 122}, 1995, No. {\bf 2},
pp. 291-308.
 
\ky{B} E. Bombieri, {\it Canonical models of general type},
Inst. Hautes Etudes Sci. Publ. Math. {\bf 42}, 1972, pp.171-219.
 
\ky{De1} J.-P. Demailly, {\it A numerical criterion for very ample
line bundles}, Journal of Differential Geometry, {\bf 37}, 1993,
pp. 323-374.
 
\ky{De2} J.-P. Demailly, { \it Effective bounds for
very ample line bundles}, Inv. Math. {\bf 124}, 1996, No. {\bf 1-3},
pp. 243-261.

\ky{De3} J.-P. Demailly, {\it $L^2$ vanishing theorems for positive line
 bundles,
and adjunction theory}, to appear.
 
\ky{EKL} L. Ein, O. K\"uchle, and
R. Lazarsfeld, {\it Local positivity of ample
line bundles}, J. Diff. Geom. Vol. {\bf 42}, No. {\bf 2}, 1995,
pp. 193-219.

\ky{EL1} L. Ein and R. Lazarsfeld, {\it Global generation of
pluricanonical and
Adjoint linear series on smooth projective threefolds},
JAMS, {\bf 6}, 1993, pp. 875-903.
 
\ky{EL2} L. Ein and R. Lazarsfeld, {\it Seshadri constants on smooth
surfaces}, Journ\'ees de G\'eometrie Alg\'ebrique
d'Orsay, Ast\'erique {\bf 218}, 1993, pp. 549-576.
 
\ky{EL3} L. Ein and R. Lazarsfeld, {\it Singulariies of theta
divisors, and the birational geometry of irregular varieties},
(to appear).
 
\ky{ELM} L.Ein, R. Lazarsfeld, and V. Masek, {\it Global generation
of linear series of terminal threeefolds}, International Journal of
Mathematics, {\bf 6}, 1995, pp. 1-18.

\ky{ELN} L. Ein, R. Lazarsfeld, and M. Nakaymaye, {\it Zero estimates,
intersection theory, and a theorem of Demailly}, Higher Dimensional
Complex varieties, Walter de Gruyter 1996.
 
\ky{EV1} H. Esnault and E. Viehweg, {\it Lectures on Vanishing
       Theorems}, DMV Seminar Band {\bf 20}, Birkh\"{a}user, 1992.
 
\ky{EV2} H. Esnault and E. Viehweg, {\it Sur une minoration du
degr\'e d'hypersurfaces s'annulant en certain points}, Math. Ann
{\bf 263}, 1983, No. {\bf 1}, pp.75-86.

\ky{Fer} G. Fern\'andez del Busto,
{\it Bogomolov's instability and Kawatmata-Viehweg's vanishing theorem},
J. of Alg. Geom., (to appear).
 
\ky{Fuj1} T. Fujita, {\it Remarks on Ein-Lazarsfeld criterion of
spannedness of adjoint bundles of polarized threefolds}, (preprint).
 
\ky{Fuj2} T. Fujita, {\it Toward a separation theorem of points
by adjoint linear systems on polarized threefolds}, (preprint).
 
\ky{Fuj3} T. Fujita, {\it Approximating Zariski decomposition of big line
bundles}, Kodai Math. J. {\bf 17} , 1994, pp. 1 - 3.
 
\ky{GL} M. Green and R. Lazarsfeld, {\it Deformation theory, generic
vanishing theorems, and some conjecture of Enriques, Catanese, and
Beauville}, Invent Math. {\bf 90} 1987, pp. 387-407.
 
\ky{H} S. Helmke, {\it On Fujita's conjecture}, (preprint).
 
\ky{Ka1} Y. Kawamata, {\it On Fujita's freeness conjecture for 3-folds
and 4-folds}, (to appear).
 
\ky{Ka2} Y. Kawamata, {\it Subadjunction of log canonical divisors
for a subvariety of codimension 2}, (to appear).
 
\ky{KMM} Y. Kawamata, K.Matsuda, and K. Matsuki, {\it Introduction
to the minimal model problems}, in Algebraic Geoemtry, Sendai Adv.
Stud. Pure Math., {\bf 10}, ed. T. Oda, Kinokuniya-North Holland,
1987, pp. 287-360.
 
\ky{KV} Y. Kawamata and E. Viehweg, {\it On a characteriazation
of abelian varieties in the classicfication theory of algebraic
varieties}, Comp. Math. {\bf 41} 1981, pp. 355-360.
 
\ky{Kol1} J. Koll\'ar, {\it Flip and abundance for algebraic
threefolds,} Ast\'erique {\bf 211}, 1992.

\ky{Kol2} J. Koll\'ar, {\it Effective basepoint freeness},
Math. Ann. {\bf 296},
1993, pp. 595-605.

\ky{Kol3} J. Koll\'ar, {\it Shafervich maps and automorphic forms},
Princeton University Press, 1995.
 
\ky{Kol4} J. Koll\'ar, {\it Singularites of pairs},
(to appear).
 
\ky{Laz} R. Lazarsfeld, {\it Lectures on linear series}, (to appear).
 
\ky{Lee} S. Lee, {\it Remarks on the pluricanonical and the
adjoint linear series}, (to appear).
 
\ky{R} I. Reider, {\it Vector bundles of rank 2 and linear systems
on surfaces}, Ann. of Math. (2) {\bf 127}, 1988, pp. 309-316.
 
\ky{Siu1} Y.-T. Siu, {\it Effective very ampleness}, Inv. Math.
{\bf 124}, 1996, No. {\bf 1-3},pp. 563-571.
 
\ky{Siu2} Y.-T. Siu, {\it Very ampleness criterion of double adjoint of ample
line bundles}, Modern Methods in Complex Analysis, Ann. of Math. Stud.,
{\bf 137}, Princeton University Press, Princeton, NJ, 1995, pp.291-318.
 
\ky{SV}  R. Smith and R. Varley,  {\it Multiplicity $g$ points on theta
divisors}, Duke Math. J. {\bf 82}, 1996, No. {\bf 2}, pp. 319-326.

\end